# Nonlinear free flexural vibrations of functionally graded rectangular and skew plates under thermal environments

N. Sundararajan[1], T. Prakash[2] and M. Ganapathi[3]*

[1]142 South car street, TV Koil, Tiruchirapalli, Tamilnadu, India.
[2] No. 11B, Sirudhunai Nayanar St, Palayamkottai, Tirunelveli Dist, Tamilnadu-627002, India.
[3] No. 925 Scientist Hostel III, DRDO Township, C.V.Raman Nagar (PO), Bangalore-560093 India.

## Abstract

The nonlinear formulation developed based on von Karman's assumptions is employed to study the free vibration of characteristics of functionally graded material (FGM) plates subjected to thermal environment. Temperature field is assumed to be a uniform distribution over the plate surface and varied in thickness direction only. Material properties are assumed to be temperature dependent and graded in the thickness direction according to simple power law distribution. The nonlinear governing equations obtained using Lagrange's equations of motion are solved using finite element procedure coupled with the direct iteration technique. The variation of nonlinear frequency ratio with amplitude is highlighted considering various parameters such as gradient index, temperature, thickness and aspect ratios, skew angle and boundary condition. For the numerical illustrations, silicon nitride/stainless steel is considered as functionally graded material. The results obtained here reveal that the temperature field and gradient index have significant effect on the nonlinear vibration of the functionally graded plate.

**Keywords**: Functionally graded plate, nonlinear free vibration, aspect ratio, temperature, gradient index.

**Running Title:** Nonlinear vibration of functionally graded plates

*Correspondence : mganapathi@rediffmail.com (M. Ganapathi)

## 1. Introduction

Functionally graded materials (FGMs) used initially as thermal barrier materials for aerospace structural applications and fusion reactors are now developed for the general use as structural components in high temperature environments and being strongly considered as a potential structural material candidate for the design of high speed aerospace vehicles. Further, these materials are inhomogeneous, in the sense that the material properties vary smoothly and continuously in one or more directions, and obtained by changing the volume fraction of the constituent materials [1-4]. For the structural integrity, FGMs are preferred over fiber-matrix composites that may result in debonding due to the mismatch in the mechanical properties across the interface of two discrete materials bonded together. Typical FGMs are made from a mixture of ceramic and metal using powder metallurgy techniques. With the increased use of these materials for structural components in many engineering applications, it is necessary to understand the dynamic characteristics of functionally graded plates.

It is seen from the literature that the amount of work carried out on the vibration characteristics of isotropic plates and composite laminates are exhaustive. However, the investigations of linear dynamic behaviors of FGM plates under thermo-mechanical environment are limited in number and are discussed briefly here. Tanigawa et al. [5] have examined transient thermal stress distribution of FGM plates induced by unsteady heat conduction with temperature dependent material properties. Reddy and Chin [6] have dealt with many problems, including transient response of plate due to heat flux. In Ref. [7], three dimensional analysis of transient thermal stress in functionally graded plates has been carried out by adopting Laplace transformation technique and power series method. He et al. [8] presented finite element formulation based on thin plate theory for the shape and vibration control of FGM plate with integrated piezoelectric sensors and actuators under mechanical load whereas Liew et al. [9] have analyzed the active vibration control of plate subjected to a thermal gradient using shear deformation theory. Ng et al. [10] have investigated the parametric resonance of plates based on Hamilton' principle and the assumed mode technique. Yang and Shen [11] have analyzed dynamic response of thin FGM plates subjected to impulsive loads using Galerkin procedure coupled with modal superposition method whereas, by neglecting the heat conduction effect, such plates and panels in thermal environments have been examined based on shear deformation with temperature dependent material properties in Ref. [12]. Cheng and Batra [13] studied the steady state vibration of a simply supported functionally graded polygonal plate with temperature independent material properties. Sills et al. [14] have presented different modeling aspects and also simulated the dynamic response under a step load.

Due to the increased utilization of thin-walled structural components in the design of flight vehicle structures, their vibration characteristics at large amplitudes in response to the conditions they are subjected to, have attracted the attention of many researchers in recent years. The studies pertaining to isotropic and composite laminates are reviewed and are well documented [15-19]. However, it is observed from the literature that such work concerning functionally graded material structures is scarce in the literature [20-22]. Praveen and Reddy [20] adopting finite element procedure analyzed the nonlinear dynamic response of functionally graded ceramic metal plates subjected to mechanical and thermal loads whereas Yang et al. [21] have investigated the large amplitude effect on free vibration characteristics of thermo-electro-mechanically stressed functionally graded plate integrated with piezoelectric actuators employing a semi analytical method involving one dimensional differential quadrature and Galerkin method. These work were based on using temperature independent material properties. Based on perturbation technique, Huang and Shen [22] studied the nonlinear free and forced vibrations of functionally graded plates in thermal environments with temperature dependent properties. It is further inferred that these available work on nonlinear dynamic of FGM plates are limited to square with amplitude-to-thickness ratio equal to 1 corresponding to the lowest vibration mode. Parametric studies including the aspect ratio, higher mode, and non-rectangular plan-form are also important while developing structural strategies with functionally graded materials.

Here, an eight-noded shear flexible quadrilateral plate element developed based on consistency approach [23, 24] is used to analyze the large amplitude free flexural vibrations behavior of FGM plates, including skew plates, subjected to thermo-mechanical environments. The geometrical non-linearity using Von Karman's assumption is introduced. The formulation also includes in-plane and rotary inertia effects. The temperature field is assumed to be constant in the plane and varied only in the thickness direction of the plate. The material is assumed to be temperature dependent and graded in the thickness direction according to the power-law distribution in terms of volume fractions of the constituents. A detailed investigation has been carried out to bring out the influences of thickness and aspect ratios and skew-angle of the plate on the frequency ratio of functionally graded plates.

## 2. Theoretical development and formulation

A functionally graded rectangular plate (length *a*, width *b,* and thickness *h*) made of a mixture of ceramics and metals is considered with the coordinates $x, y$ along the in-pane directions and $z$ along the thickness direction. The material in top surface $(z = h/2)$ of the plate and in bottom surface $(z = -h/2)$ of the plate is ceramic and metal, respectively. The effective material properties $P$, such as Young's modulus $E$, and thermal expansion coefficient $\alpha$, can be written as [20]

$$P = P_c V_c + P_m V_m \tag{1}$$

where $P_c$ and $P_m$ are the material properties of the ceramic rich top surface and metal rich bottom surface, respectively. $V_c$ and $V_m$ are volume-fractions of ceramic and metal respectively and are related by

$$V_c + V_m = 1 \quad (2)$$

The properties of the plate are assumed to vary through the thickness. The property variation is assumed to be in terms of a simple power law. The volume fraction $V_c$ is expressed as

$$V_c(z) = \left(\frac{2z+h}{2h}\right)^k \quad (3)$$

where $k$ is the volume fraction exponent $(k \geq 0)$. The material properties *P* that are temperature dependent can be written as

$$P = P_0(P_{-1}T^{-1} + 1 + P_1 T + P_2 T^2 + P_3 T^3) \quad (4)$$

where $P_0$, $P_{-1}$, $P_1$, $P_2$ and $P_3$ are the coefficients of temperature $T(K)$ and are unique to each constituent.

From Eqs. (1) - (4), the modulus of elasticity $E$, the coefficient of thermal expansion $\alpha$, the density $\rho$ and the thermal conductivity $K$ are written as

$$\begin{aligned}
E(z,T) &= (E_c(T) - E_m(T))\left(\frac{2z+h}{2h}\right)^k + E_m(T) \\
\alpha(z,T) &= (\alpha_c(T) - \alpha_m(T))\left(\frac{2z+h}{2h}\right)^k + \alpha_m(T) \\
\rho(z) &= (\rho_c - \rho_m)\left(\frac{2z+h}{2h}\right)^k + \rho_m \\
K(z) &= (K_c - K_m)\left(\frac{2z+h}{2h}\right)^k + K_m
\end{aligned} \quad (5)$$

Here the mass density $\rho$ and thermal conductivity $K$ are assumed to be independent of temperature. The Poisson's ratio $\nu$ is assumed to be a constant $\nu(z) = \nu_0$.

The temperature variation is assumed to occur in the thickness direction only and the temperature field is considered constant in the *xy* plane. In this case, the temperature through thickness is governed by the one-dimensional Fourier equation of heat conduction:

$$\frac{d}{dz}\left[K(z)\frac{dT}{dz}\right]=0, \qquad T=T_c \ at \ z=h/2$$
$$T=T_m \ at \ z=-h/2 \tag{6}$$

The solution of Eq. (6) is obtained by means of polynomial series [25] and given by

$$T(z) = T_m + (T_c - T_m)\eta(z) \tag{7}$$

$$where \quad \eta(z)=\frac{1}{C}\left[\left(\frac{2z+h}{2h}\right)-\frac{K_{cm}}{(k+1)K_m}\left(\frac{2z+h}{2h}\right)^{k+1}+\frac{K_{cm}^2}{(2k+1)K_m^2}\left(\frac{2z+h}{2h}\right)^{2k+1}-\right.$$
$$\left.\frac{K_{cm}^3}{(3k+1)K_m^3}\left(\frac{2z+h}{2h}\right)^{3k+1}+\frac{K_{cm}^4}{(4k+1)K_m^4}\left(\frac{2z+h}{2h}\right)^{4k+1}-\frac{K_{cm}^5}{(5k+1)K_m^5}\left(\frac{2z+h}{2h}\right)^{5k+1}\right];$$

$$C=1-\frac{K_{cm}}{(k+1)K_m}+\frac{K_{cm}^2}{(2k+1)K_m^2}-\frac{K_{cm}^3}{(3k+1)K_m^3}+\frac{K_{cm}^4}{(4k+1)K_m^4}-\frac{K_{cm}^5}{(5k+1)K_m^5};$$

and $K_{cm} = K_c - K_m$

Using Mindlin formulation, the displacements $u, v, w$ at a point $(x, y, z)$ in the plate (Fig. 1a) from the medium surface are expressed as functions of mid-plane displacements $u_0, v_0$ and $w$, and independent rotations $\theta_x$ and $\theta_y$ of the normal in $xz$ and $yz$ planes, respectively, as

$$\begin{aligned} u(x,y,t) &= u_0(x,y,t) + z\theta_x(x,y,t) \\ v(x,y,t) &= v_0(x,y,t) + z\theta_y(x,y,t) \\ w(x,y,t) &= w_0(x,y,t) \end{aligned} \tag{8}$$

where t is the time.

von Karman's assumptions for moderately large deformation allows Green's strains to be written in terms of mid-plane deformation of Eqn. (8) for a plate as,

$$\{\varepsilon\}=\{\varepsilon^L\}+\{\varepsilon^{NL}\} \tag{9}$$

Taking into account the effect of shear deformation, the total linear and nonlinear strain at any point can be expressed as

$$\{\varepsilon^{L}\} = \begin{Bmatrix} \varepsilon_p^{L} \\ 0 \end{Bmatrix} + \begin{Bmatrix} z\,\varepsilon_b \\ \varepsilon_s \end{Bmatrix}$$ and

$$\{\varepsilon^{NL}\} = \begin{Bmatrix} \varepsilon_p^{NL} \\ 0 \end{Bmatrix} \tag{10}$$

The mid-plane strains $\{\varepsilon_p^{L}\}$, bending strains $\{\varepsilon_b\}$, shear strains $\{\varepsilon_s\}$ and the nonlinear components of in-plane strains $\{\varepsilon_p^{NL}\}$ in Eqn. (10) are written as

$$\{\varepsilon_p^{L}\} = \begin{Bmatrix} u_{o,x} \\ v_{o,y} \\ u_{o,y} + v_{o,x} \end{Bmatrix} \tag{11}$$

$$\{\varepsilon_b\} = \begin{Bmatrix} \theta_{x,x} \\ \theta_{y,y} \\ \theta_{x,y} + \theta_{y,x} \end{Bmatrix} \tag{12}$$

$$\{\varepsilon_s\} = \begin{Bmatrix} \theta_x + w_{o,x} \\ \theta_y + w_{o,y} \end{Bmatrix} \tag{13}$$

$$\{\varepsilon_p^{NL}\} = \begin{Bmatrix} (1/2)w_{o,x}^2 \\ (1/2)w_{o,x}^2 \\ w_{o,x}\,w_{o,y} \end{Bmatrix} \tag{14}$$

where the subscript comma denotes the partial derivative with respect to the spatial coordinate succeeding it.

The membrane stress resultants $\{N\}$ and the bending stress resultants $\{M\}$ can be related to the membrane strains $\{\varepsilon_p\}$ ($= \{\varepsilon_p^{L}\} + \{\varepsilon_p^{NL}\}$) and bending strains $\{\varepsilon_b\}$ through the constitutive relations by

$$\{N\} = \begin{Bmatrix} N_{xx} \\ N_{yy} \\ N_{xy} \end{Bmatrix} = [A_{ij}]\{\varepsilon_p\} + [B_{ij}]\{\varepsilon_b\} - \{N^T\} \tag{15}$$

$$\{M\} = \begin{Bmatrix} M_{xx} \\ M_{yy} \\ M_{xy} \end{Bmatrix} = [B_{ij}]\{\varepsilon_p\} + [D_{ij}]\{\varepsilon_b\} - \{M^T\} \tag{16}$$

where the matrices $[A_{ij}]$ ,$[B_{ij}]$ and $[D_{ij}]$ $(i, j = 1,2,6)$ are the extensional, bending-extensional coupling and bending stiffness coefficients and are defined as $[A_{ij}, B_{ij}, D_{ij}] = \int_{-h/2}^{h/2} [\overline{Q}_{ij}](1, z, z^2)dz$ . The thermal stress resultant $\{N^T\}$ and moment resultant $\{M^T\}$ are

$$\{N^T\} = \begin{Bmatrix} N_x^T \\ N_y^T \\ N_{zy}^T \end{Bmatrix} = \int_{-h/2}^{h/2} [\overline{Q}_{ij}] \begin{Bmatrix} \alpha(z,T) \\ \alpha(z,T) \\ 0 \end{Bmatrix} \Delta T(z)\, dz \tag{17}$$

$$\{M^T\} = \begin{Bmatrix} M_x^T \\ M_y^T \\ M_{zy}^T \end{Bmatrix} = \int_{-h/2}^{h/2} [\overline{Q}_{ij}] \begin{Bmatrix} \alpha(z,T) \\ \alpha(z,T) \\ 0 \end{Bmatrix} \Delta T(z)\, zdz \tag{18}$$

where the thermal coefficient of expansion $\alpha(z,T)$ is given by Eq. (5), and $\Delta T(z) = T(z) - T_0$ is temperature rise from the reference temperature $T_0$ at which there are no thermal strains.

Similarly the transverse shear force $\{Q\}$ representing the quantities $\{Q_{xz}, Q_{yz}\}$ is related to the transverse shear strains $\{\varepsilon_s\}$ through the constitutive relations as

$$\{Q\} = [E_{ij}]\{\varepsilon_s\} \tag{19}$$

where $$E_{ij} = \int_{-h/2}^{h/2} [\overline{Q}_{ij}]\kappa_i \kappa_j dz$$

Here $[E_{ij}]$ $(i, j = 4,5)$ are the transverse shear stiffness coefficients, $\kappa_i$ is the transverse shear coefficient for non-uniform shear strain distribution through the plate thickness. $\overline{Q}_{ij}$ are the stiffness coefficients and are defined as

$$\overline{Q}_{11} = \overline{Q}_{22} = \frac{E(z,T)}{1-\nu^2}; \ \overline{Q}_{12} = \frac{\nu E(z,T)}{1-\nu^2}; \ \overline{Q}_{16} = \overline{Q}_{26} = 0; \ \overline{Q}_{44} = \overline{Q}_{55} = \overline{Q}_{66} = \frac{E(z,T)}{2(1+\nu)}$$

(20)

where the modulus of elasticity $E(z,T)$ is given by Eq. (5).

The strain energy functional $U$ is given as

$$U(\delta) = (1/2)\int_A \Big[\{\varepsilon_p\}^T [A_{ij}]\{\varepsilon_p\} + \{\varepsilon_p\}^T [B_{ij}]\{\varepsilon_b\} + \{\varepsilon_b\}^T [B_{ij}]\{\varepsilon_p\} + \{\varepsilon_b\}^T [D_{ij}]\{\varepsilon_b\} + \{\varepsilon_s\}^T [E_{ij}]\{\varepsilon_s\} - \{\varepsilon_p^0\}^T \{N^T\} - \{\varepsilon_b\}^T \{M^T\}\Big] dA$$

(21)

where $\delta$ is the vector of the degree of freedom associated to the displacement filed in a finite element discretisation.

Following the procedure given in Ref. [26], the potential energy functional U given in Eqn. (21) can be rewritten as

$$U(\delta) = \{\delta\}^T [(1/2)[\mathbf{K}] + [(1/6)[\mathbf{N_1}(\delta)] + (1/12)[\mathbf{N_2}(\delta)] + (1/2)[\mathbf{N_3}]]\{\delta\}, \quad (22)$$

where [**K**] is the linear stiffness matrix of the laminate. [$\mathbf{N_1}$] and [$\mathbf{N_2}$] are non-linear stiffness matrices linearly and quadratically dependent on the field variables, respectively. [$\mathbf{N_3}$] is transverse shear stiffness matrix of the plate.

The kinetic energy of the plate is given by

$$T(\delta) = (1/2)\int_A \left[p\left(\dot{u}_0^2 + \dot{v}_0^2 + \dot{w}_0^2\right) + I\left(\dot{\theta}_x^2 + \dot{\theta}_y^2\right)\right] dA \quad (23)$$

where $p = \int_{-h/2}^{h/2} \rho(z)dz$, $I = \int_{-h/2}^{h/2} z^2 \rho(z)dz$ and $\rho(z)$ is mass density which varies through the thickness of the plate and is given by Eq. (5).

The plate is subjected to temperature filed and this, in turn, results in-plane stress resultants $\left(N_{xx}^{th}, N_{yy}^{th}, N_{xy}^{th}\right)$. Thus, the potential energy due to pre-buckling stresses $\left(N_{xx}^{th}, N_{yy}^{th}, N_{xy}^{th}\right)$ developed under thermal load can be written as

$$V(\delta) = \int_A \left\{ \frac{1}{2}\left[ N_{xx}^{th}\left(\frac{\partial w}{\partial x}\right)^2 + N_{yy}^{th}\left(\frac{\partial w}{\partial y}\right)^2 + 2N_{xy}^{th}\left(\frac{\partial w}{\partial x}\right)\left(\frac{\partial w}{\partial y}\right)\right] + \frac{h^3}{24}\left[ N_{xx}^{th}\left\{\left(\frac{\partial \theta_x}{\partial x}\right)^2 + \left(\frac{\partial \theta_y}{\partial x}\right)^2\right\} + N_{yy}^{th}\left\{\left(\frac{\partial \theta_x}{\partial y}\right)^2 + \left(\frac{\partial \theta_y}{\partial y}\right)^2\right\} + 2N_{xy}^{th}\left\{\left(\frac{\partial \theta_x}{\partial x}\right)\left(\frac{\partial \theta_x}{\partial y}\right) + \left(\frac{\partial \theta_y}{\partial x}\right)\left(\frac{\partial \theta_y}{\partial y}\right)\right\}\right] \right\} dA$$

(24)

Substituting Eqs. (22-24) in Lagrange's equation of motion, one obtains the governing equations as

$$[\mathbf{M}]\{\ddot{\delta}\} + ([\mathbf{K}] + (1/2)[N_1] + (1/3)[N_2] + [N_3] + [\mathbf{K_G}])\{\delta\} = \{\mathbf{o}\} \quad (25)$$

where **[M]** is the consistent mass matrix; **[K$_G$]** is the geometric stiffness matrix, respectively. **{ $\ddot{\delta}$ }** is the acceleration vector.

Substituting characteristics of the time function at the point of reversal of motion [24]

$$\{\ddot{\delta}\} = -\omega^2\{\delta\} \quad (26)$$

Eq. (25) will lead to the following nonlinear algebraic equation form,

$$([\mathbf{K}] + (1/2)[N_1] + (1/3)[N_2] + [N_3] + [\mathbf{K_G}])\{\delta\} - \omega^2[M]\{\delta\} = \{\mathbf{o}\} \quad (27)$$

where, $\omega$ is the natural frequency. The frequency-amplitude relation is obtained by solving Eq. (27) through finite element procedure in conjunction with direct iteration technique.

## 3. Element description

The plate element employed here is a $C^0$ continuous shear flexible element and needs five nodal degrees of freedom $u_0, v_0, w_o, \theta_x, \theta_y$ at eight nodes in QUAD-8 element. If the interpolation functions for QUAD-8 are used directly to interpolate the five variables $u_0$ to $\theta_y$ in deriving the shear strains and membrane strains, the element will lock and show oscillations in the shear and membrane stresses. Field consistency requires that the transverse shear strains and membrane strains must be interpolated in a consistent manner. Thus $\theta_x$ and $\theta_y$ terms in the expressions for $\{\varepsilon_s\}$ given by Eq. (13) have to be consistent with field functions $w_{o,x}$ and $w_{o,y}$. This is achieved by using field redistributed substitute shape functions to interpolate those specific terms, which must be consistent, as described in Ref [23, 24]. This element is free from locking syndrome and has good convergence properties. For the sake of brevity, these are not presented here, as they are available in the literature [23, 24]. Since the element is based on field consistency approach, exact integration is applied for calculating various strain energy terms.

## 4. Skew Boundary Transformation

For skew plates supported on two adjacent edges, the edges of the boundary elements may not be parallel to the global axes $(x, y, z)$. In such a situation, it is not possible to specify the boundary conditions in terms of the global displacements $u_o, v_o, w_o$, etc. In order to specify the boundary conditions at skew edges, it is necessary to use edge displacements $u_o^l, v_o^l, w_o^l$, etc. in local coordinates $(x^l, y^l, z^l)$ as shown in Fig. 1b. It is thus required to transform the element matrices corresponding to global axes to local edge axes with respect to which the boundary conditions can be conveniently specified. The relation between the global and local degrees of freedom of a node can be obtained through the simple transformation rules [27] and the same can be expressed as

$$d_i = L_g d_i^l \tag{28.a}$$

in which $d_i, d_i^l$ are generalized displacement vectors in the global and local coordinate system, respectively of node $i$ and they are defined as

$$d_i = [u_o \quad v_o \quad w_o \quad \theta_x \quad \theta_y]^T \tag{28.b}$$

$$d_i^l = [u_o^l \quad v_o^l \quad w_o^l \quad \theta_x^l \quad \theta_y^l]^T \tag{28.c}$$

The nodal transformation matrix for a node *i*, on the skew boundary is

$$L_g = \begin{bmatrix} c & s & 0 & 0 & 0 \\ -s & c & 0 & 0 & 0 \\ 0 & 0 & 1 & 0 & 0 \\ 0 & 0 & 0 & c & s \\ 0 & 0 & 0 & -s & c \end{bmatrix} \tag{29}$$

in which $c = \cos(\psi)$ and $s = \sin(\psi)$, where $\psi$ is the angle of the plate. It may be noted that for the nodes, which are not lying on the skew edges, the node transformation matrix has only nonzero values for the principal diagonal elements, which are equal to 1. Thus, for the complete element, the element transformation matrix is written as

$$[T]_e = diag \left\langle L_g \quad L_g \quad L_g \quad L_g \quad L_g \quad L_g \quad L_g \quad L_g \right\rangle \tag{30}$$

For those elements whose nodes are on the skew edges, the element matrices are transformed to the local axes using the element transformation matrix $T_e$ and then the global matrices/vectors are obtained using standard assembly procedures.

## 5. Solution Procedure

The vibration problem is solved using eigenvalue formulation. To solve the non-linear eigenvalue problems, an iterative procedure is used. Firstly, the eigenvector (mode shape) is obtained from the linear vibration analysis, neglecting the non-linear stiffness matrix in Eq.(27) and then normalized. Next, the normalized vector is amplified/scaled up so that the maximum displacement is equal to the desired amplitude, say $w/h = 0.2$ ( $w$ ) is the maximum lateral displacement, $h$ is the thickness of the plate). This gives the initial vector, denoted by $\overline{\delta}$ . The iterative solution procedure for the non-linear analysis starts with this initial vector. Based on this initial mode shape ($\overline{\delta}$ ), the non-linear stiffness matrix that depends on displacement (linearly and quadratically) is formed. Subsequently, the updated eigenvalue and its corresponding eigenvector are obtained. This eigenvector is further normalized, and scaled up by the same amplitude $(w/h)$, and the iterative procedure adopted here continues till the frequency values and mode shapes evaluated from the subsequent two iterations satisfy the prescribed convergence criteria [28] as

$$\sum_{m}^{N} \frac{(\omega_m^r - \omega_m^{r-1})}{\omega_m^r} \leq 0.0001$$

$$\sum_{i}^{N} \frac{(\delta_i^r - \delta_i^{r-1})}{\delta_i^r} \leq 0.0001 \text{, for } m^{\text{th}} \text{ mode,} \tag{31}$$

where, $m$ , $i$, $N$ and $r$ represents the mode number, degree of freedom of the finite element model, total degree of freedom and iteration number, respectively.

## 6. Results and Discussion

The study, here, has been focused on the large amplitude free flexural vibration behavior of functionally graded plates. Fig. 2 shows the variation of the volume fractions of ceramic and metal respectively, in the thickness direction *z* for the FGM plate. The top surface is ceramic rich and the bottom surface is metal rich. The FGM plate considered here consists of Silicon nitride ($Si_3N_4$) and stainless steel (SUS304). The temperature coefficients corresponding to Si3N4 / SUS304 are listed in Table 1 [6]. The mass density and thermal conductivity are: $\rho_c = 2370\ kg/m^3$, $K_c = 9.19\ W/mK$ for $Si_3N_4$; and $\rho_m = 8166\ kg/m^3$, $K_m = 12.04\ W/mK$ for SUS304. Poisson's ratio $\nu$ is assumed to be a constant and equals to 0.28. Transverse shear coefficient is taken as 0.91. The plate is of uniform thickness and boundary conditions considered here are:

simply supported :

$u_o = w_o = \theta_y = 0$ on $x = 0, a$ and $v = w = \theta_x = 0$ on $y = 0, b$

clamped support :

$u = v = w = \theta_x = \theta_y = 0$ on $x = 0, a$ & $y = 0, b$

Based on the progressive mesh refinement, an 8x8 mesh is found to be adequate to model the full plate for the present analysis [29]. Before proceeding for the detailed study for the large amplitude free flexural vibrations behavior of functionally graded plates, the formulation developed herein is validated against the available results [22,30,31] pertaining to the linear case of FGM plates and nonlinear vibrations of isotropic plates in Tables 2a and 2b, respectively. Here, the calculated non-dimensional linear frequency is defined as: $\omega = \omega_L \left(\frac{a^2}{h}\right)\left(\frac{\rho_m (1-\nu^2)}{E_m}\right)^{1/2}$, where $\rho_m$ and $E_m$ are the mass density and Young's modulus of metal, respectively. The results are found to be in good agreement with the existing solutions. It is observed from Table 2a that, with the increase in power law index *k* up to certain value, the rate of decrease in the frequency value is high, and further increase in *k* leads to less reduction in the frequency. For the low values of *k*, the stiffness degradation occurs due to the increase in the metallic volumetric fraction.

Next, the nonlinear flexural vibration behavior of FGM is numerically studied with and without thermal environment. For the uniform temperature case, the material properties are evaluated at $T = 300\,K$. The variation of nonlinear to linear frequency ratio ($\omega_{NL}/\omega_L$, where subscripts NL and L corresponds to the non-linear and linear respectively) with respect to non-dimensional amplitude ($w/h$; $w$ is the flexural amplitude of the plate) evaluated for simply supported plate is shown in Table 3. It is seen from Table 3 that the frequency ratios decreases with the increase in the gradient index *k* up to certain value, i. e. *k*=2, and then the frequency ratio increases with further increase in *k*. This observation is true, irrespective of thick or thin, and square or rectangular case considered in the present study. For the higher values of *k*, it is revealed from the detailed analysis that, although the linear as well as nonlinear frequency values decrease, the reduction in linear frequency is more compared to that of the nonlinear case and, for the sake of brevity, it is not shown here. Due to greater reduction in linear frequency values at higher *k*, the overall trend in nonlinear frequency ratio thus increases. It can be further viewed that the frequency ratio increases with the increase in amplitude, as expected. However, it is noticed that there is a sudden drop in the increasing frequency trend at certain higher amplitude and then gradually increases with further increase in amplitude exhibiting hardening behavior. This is possibly attributed to the change in stiffness values, and thus leading to the redistribution of mode shapes associated with certain level of amplitudes of vibration, as highlighted in Fig. 3 loosing symmetry and shifting the maximum displacement towards one side of the plate. It may be opined from Table 3 that, in general, this trend of dropping in frequency occurs around $w/h$ = 1.2 or higher for the square plate whereas it corresponds to around $w/h$ = 0.8 or less for rectangular plate considered here. The abrupt drop in the frequency ratio with respect to gradient index and amplitude ratio is highlighted in the Table. Also it can be concluded that, for thin case, the reduction in the frequency ratio value occurs at fairly high amplitude compared to thick case. Furthermore, it may be inferred that the frequency ratio decreases with the increase in the thickness ratio (a/h), and there is qualitatively no change in the variation of frequency, irrespective of thickness parameter.

The influence of thermal gradient on the nonlinear vibrations characteristics of FGMs is examined in Tables 4 & 5 with different surface temperatures. The temperature field is assumed to vary only in the thickness direction and determined by the expression given in Eq. (7). The temperature for the ceramic surface is varied ($T_c$ = 400*K*, 600*K*) while keeping the constant value for metallic surface ($T_m$ = 300*K*). To evaluate the non-dimensional frequency given in Eq. (25), $\rho_m$ and $E_m$ are taken at $T_0 = 300\ K$ . The frequency-drops are also highlighted in these Tables. The variation of linear and nonlinear frequencies is qualitatively similar to those of non-thermal case. However, the nonlinear frequency ratios are higher here compared to the non-thermal case given in Table 3, irrespective of thickness and aspect ratios of the FGM plates. Further more, it is observed from Tables 4 and 5 that the increase in temperature gradient results in relatively higher nonlinear frequency ratios.

Similar study is carried out for clamped FGM plates subjected to the thermal gradient ($T_c = 600\ K$ and $T_m = 300\ K$) and the nonlinear frequency ratios obtained for square and rectangular cases are tabulated in Table 6. It is evident from this Table that, in comparison to simply supported case, the nonlinear frequency ratios are, in general, less. It is further viewed that the lowest frequency value occurs corresponding to higher gradient index of *k*=5, irrespective of aspect ratio, unlike *k*=2 for simply supported condition.

Lastly, the effect of simply supported skewed FGM plates is investigated and the results are shown in Table 7 for various amplitude and gradient index. The general nonlinear behavior is qualitatively same as that of square plate. However, it is observed that the nonlinear frequency ratios are, in general, higher compared to the rectangular case. It is also revealed that with the increase in the skew angle, the drop in frequency ratio occurs at relatively low amplitude that is associated with mode redistribution.

## 7. Conclusions

The nonlinear flexural vibration behavior of FGM plates under thermal environment has been studied. The formulation is based on first-order shear deformation theory and includes geometric nonlinear using von Karman's assumptions. The material properties are assumed to be varied through the thickness direction based on power law distribution and the temperature varies only in the thickness direction. The numerical experiments have been conducted to bring out the effectiveness of gradient index, aspect and thickness ratios, boundary conditions and thermal environment on the nonlinear flexural vibrations FGM plates. From the detailed parametric study, the following observations can be made:

1. With the increase in gradient index, *k*, the fundamental frequency decreases due to degradation of stiffness by the metallic inclusion.
2. With the increase in gradient index *k* value, the nonlinear frequency ratio decreases initially and then increases.
3. The increase in the nonlinear frequency ratio drops at certain amplitude (*w/h*) due to mode redistribution and then increases with further increase in amplitude.
4. The sudden drop in frequency occurs at low amplitude with the increase in the aspect ratio whereas it happens at fairly high amplitude for thin case.
5. Although the thermal environment reduces the natural frequency value as expected, the nonlinear frequency ratio is high compared to non-thermal case. Also higher temperature gradient yields, in general, higher nonlinear frequency ratios.
6. Skew effect, in general, enhances the nonlinear frequency ratios compared to rectangular case.
7. With the increase in the skew angle, the drop in frequency ratio occurs at relatively low amplitude.

**Legends for Tables**

Table1. Temperature dependent coefficients for material Si3N4 / SUS304. Ref. [6]

Table 2a. Comparison of non-dimensional linear frequencies of simply supported FGM plate (a/b =1, a/h=8)

Table 2b. Comparison of nonlinaer frequency ratios of isotropic simply supported square plate (a/h=1000)

Table 3. Effect of aspect and thickness ratios on the nonlinear frequency ratio ($\omega_{NL}/\omega_L$) of simply supported FGM plates in ambient temperature ($T_c$=300K, $T_m$=300K)

Table 4. Effect of aspect and thickness ratios on the nonlinear frequency ratio ($\omega_{NL}/\omega_L$) of Simply supported FGM plates with temperature gradient ($T_c$ =400K, $T_m$ =300K)

Table 5. Eeffect of aspect and thickness ratios on the nonlinear frequency ratio ($\omega_{NL}/\omega_L$) of Simply supported FGM plates with temperature gradient ($T_c$ =600K, $T_m$ =300K)

Table 6. Influence of aspect ratio on the nonlinear frequency ratio ($\omega_{NL}/\omega_L$) of clamped FGM plates (a/h=20) with temperature gradient ($T_c$ =400K, $T_m$ =300K)

Table 7. Influence of skew angle of simply supported FGM plate (a/h=10, a/b=1) with temperature gradient ($T_c$ =400K, $T_m$21=300K) on the nonlinear frequency ratio ($\omega_{NL}/\omega_L$)

## Legends for Figures

Fig 1a. Configuration and coordinate system of a rectangular FGM plate.

1b. Coordinate system of a skew plate.

Fig 2. Variation of volume fractions through thickness: a) Ceramic; b) Metal.

Fig 3. The redistribution of normalized nonlinear mode shape contours of simply supported FGM plate (a/b=1, a/h=10, k=2): (a) w/h= 1.0; (b) w/h= 1.4; (c) mode shapes along y = b/2; (d) mode shapes along x = a/2.

**Table1. Temperature dependent coefficients for material $Si_3N_4$ / SUS304, Ref. [6].**

| Materials | Properties | $P_0$ | $P_{-1}$ | $P_1$ | $P_2$ | $P_3$ | P (T=300K) |
|---|---|---|---|---|---|---|---|
| $Si_3N_4$ | E (Pa) | 348.43e+9 | 0.0 | -3.070e-4 | 2.160e-7 | -8.946e-11 | 322.2715e+9 |
| | $\alpha$ (1/K) | 5.8723e-6 | 0.0 | 9.095e-4 | 0.0 | 0.0 | 7.4746e-6 |
| SUS304 | E (Pa) | 201.04e+9 | 0.0 | 3.079e-4 | -6.534e-7 | 0.0 | 207.7877e+9 |
| | $\alpha$ (1/K) | 12.330e-6 | 0.0 | 8.086e-4 | 0.0 | 0.0 | 15.321e-6 |

**Table 2a. Comparison of non-dimensional linear frequencies of simply supported FGM plate (a/b =1, a/h=8)**

| Temperature | k | Mode (1,1) | | Mode (1,2) | | Mode (2,2) | |
|---|---|---|---|---|---|---|---|
| | | Ref [22] | Present | Ref [22] | Present | Ref [22] | Present |
| Tc=400 Tm=300 | 0.0 | 12.397 | 12.311 | 29.083 | 29.016 | 43.835 | 44.094 |
| | 0.5 | 8.615 | 8.483 | 20.215 | 19.979 | 30.530 | 30.391 |
| | 1.0 | 7.474 | 7.444 | 17.607 | 17.511 | 26.590 | 26.648 |
| | 2.0 | 6.693 | 6.679 | 15.762 | 15.706 | 23.786 | 23.894 |
| | 10.0 | --- | 5.742 | --- | 13.560 | --- | 20.609 |
| Tc=600 Tm=300 | 0.0 | 11.984 | 11.888 | 28.504 | 28.421 | 43.107 | 43.343 |
| | 0.5 | 8.269 | 8.150 | 19.784 | 19.534 | 29.998 | 29.836 |
| | 1.0 | 7.171 | 7.131 | 17.213 | 17.101 | 26.109 | 26.139 |
| | 2.0 | 6.398 | 6.376 | 15.384 | 15.314 | 23.327 | 23.410 |
| | 10.0 | --- | 5.423 | --- | 13.146 | --- | 20.100 |

**Table 2b. Comparison of nonlinaer frequency ratios of isotropic simply supported square plate (a/h=1000)**

| w/h | Analytical | FEM | |
|---|---|---|---|
| | Ref [30] | Ref [31] | Present |
| 0.2 | 1.02599 | 1.02542 | 1.02563 |
| 0.4 | 1.10027 | 1.09821 | 1.09918 |
| 0.6 | 1.21402 | 1.21190 | 1.21258 |
| 0.8 | 1.35735 | 1.35588 | 1.35659 |
| 1.0 | 1.52192 | 1.52005 | 1.52339 |

**Table 3. Effect of aspect and thickness ratios on the nonlinear frequency ratio ($\omega_{NL}/\omega_L$) of simply supported FGM plates in ambient temperature ($T_c$=300K, $T_m$=300K).**

| a/b | a/h | k | *w/h* | | | | | | | |
|---|---|---|---|---|---|---|---|---|---|---|
| | | | 0.2 | 0.4 | 0.6 | 0.8 | 1.0 | 1.2 | 1.4 | 1.6 |
| 1 | 10 | 0 | 1.0271 | 1.1047 | 1.2240 | 1.3749 | 1.5490 | **1.7397*** | 1.6206 | 1.6656 |
| | | 0.5 | 1.0120 | 1.0774 | 1.1882 | 1.3344 | 1.5062 | **1.6965** | 1.5746 | 1.6143 |
| | | 1 | 1.0066 | 1.0664 | 1.1724 | 1.3144 | 1.4830 | **1.6703** | 1.5582 | 1.6057 |
| | | *2* | *1.0048* | *1.0617* | *1.1638* | *1.3017* | *1.4661* | ***1.6496*** | *1.5458* | *1.5871* |
| | | 5 | 1.0095 | 1.0691 | 1.1721 | 1.3090 | 1.4714 | **1.6520** | 1.5588 | 1.5998 |
| | | 10 | 1.0151 | 1.0798 | 1.1868 | 1.3269 | 1.4913 | **1.6731** | 1.5820 | 1.6262 |
| 1 | 20 | 0 | 1.0261 | 1.1009 | 1.2161 | 1.3622 | 1.5311 | 1.7165 | **1.9119** | 1.6966 |
| | | 0.5 | 1.0110 | 1.0736 | 1.1803 | 1.3216 | 1.4883 | **1.6732** | 1.6042 | 1.6509 |
| | | 1 | 1.0056 | 1.0626 | 1.1644 | 1.3016 | 1.4650 | **1.6472** | 1.5833 | 1.6246 |
| | | *2* | *1.0039* | *1.0579* | *1.1558* | *1.2888* | *1.4480* | *1.6265* | ***1.8154*** | *1.6141* |
| | | 5 | 1.0085 | 1.0653 | 1.1641 | 1.2961 | 1.4531 | 1.6282 | **1.8151** | 1.6319 |
| | | 10 | 1.0141 | 1.0760 | 1.1788 | 1.3139 | 1.4731 | 1.6496 | **1.8391** | 1.6767 |
| 1 | 100 | 0 | 1.0258 | 1.0996 | 1.2136 | 1.3582 | 1.5256 | 1.7091 | **1.9056** | 1.7045 |
| | | 0.5 | 1.0107 | 1.0723 | 1.1778 | 1.3175 | 1.4827 | 1.6658 | **1.8619** | 1.6584 |
| | | 1 | 1.0053 | 1.0614 | 1.1619 | 1.2973 | 1.4593 | 1.6403 | **1.8342** | 1.6356 |
| | | 2 | 1.0036 | 1.0566 | 1.1533 | 1.2844 | 1.4422 | 1.6189 | **1.8066** | 1.6278 |
| | | 5 | 1.0082 | 1.0641 | 1.1615 | 1.2918 | 1.4470 | 1.6206 | **1.8075** | 1.6393 |
| | | 10 | 1.0138 | 1.0747 | 1.1763 | 1.3096 | 1.4673 | 1.6420 | **1.8300** | 1.6593 |
| 2 | 10 | 0 | 1.0363 | 1.1398 | 1.2993 | **1.4998** | 1.3159 | 1.3482 | 1.3794 | 1.4060 |
| | | 0.5 | 1.0176 | 1.1063 | 1.2553 | **1.4508** | 1.2850 | 1.3173 | 1.3476 | 1.3735 |
| | | 1 | 1.0109 | 1.0928 | 1.2359 | **1.4260** | 1.2721 | 1.3038 | 1.3325 | 1.3588 |
| | | *2* | *1.0087* | *1.0870* | *1.2254* | ***1.4115*** | *1.2688* | *1.2972* | *1.3271* | *1.3521* |
| | | 5 | 1.0146 | 1.0962 | 1.2355 | **1.4203** | 1.2774 | 1.3063 | 1.3362 | 1.3613 |
| | | 10 | 1.0214 | 1.1092 | 1.2535 | **1.4415** | 1.2885 | 1.3198 | 1.3499 | --- |
| 2 | 20 | 0 | 1.0331 | 1.1279 | 1.2739 | **1.4596** | 1.3312 | 1.3632 | 1.3971 | |
| | | 0.5 | 1.0145 | 1.0944 | 1.2305 | **1.4107** | 1.2997 | 1.3324 | 1.3668 | 1.3992 |
| | | 1 | 1.0079 | 1.0810 | 1.2111 | **1.3864** | 1.2866 | 1.3191 | 1.3534 | 1.3862 |
| | | *2* | *1.0057* | *1.0751* | *1.2004* | ***1.3705*** | *1.2812* | *1.3130* | *1.3467* | *1.3796* |
| | | 5 | 1.0114 | 1.0842 | 1.2103 | **1.3791** | 1.2912 | 1.3222 | 1.3549 | 1.3882 |
| | | 10 | 1.0183 | 1.0972 | 1.2285 | 1.4008 | 1.5953 | 1.8778 | 1.3681 | 1.7600 |
| 2 | 100 | 0 | 1.0321 | 1.1242 | 1.2658 | **1.4468** | 1.3382 | 1.3701 | 1.4041 | 1.4383 |
| | | 0.5 | 1.0135 | 1.0908 | 1.2226 | **1.3978** | 1.3064 | 1.3391 | 1.3741 | 1.4092 |
| | | 1 | 1.0069 | 1.0773 | 1.2033 | **1.3730** | 1.2931 | 1.3258 | 1.3607 | 1.3961 |
| | | *2* | *1.0047* | *1.0714* | *1.1926* | ***1.3571*** | *1.2878* | *1.3197* | *1.3548* | *1.3886* |
| | | 5 | 1.0099 | 1.0805 | 1.2024 | **1.3656** | 1.2981 | 1.3291 | 1.3621 | 1.3967 |
| | | 10 | 1.0173 | 1.0936 | 1.2204 | **1.3872** | 1.3113 | 1.3423 | 1.3789 | 1.4094 |

*Denotes change in trend

**Table 4. Effect of aspect and thickness ratios on the nonlinear frequency ratio ($\omega_{NL}/\omega_L$) of simply supported FGM plates with temperature gradient ($T_c$ =400K, $T_m$ =300K)**

| a/b | a/h | k | *w/h* 0.2 | 0.4 | 0.6 | 0.8 | 1.0 | 1.2 | 1.4 | 1.6 |
|---|---|---|---|---|---|---|---|---|---|---|
| 1 | 10 | 0 | 1.0290 | 1.1099 | 1.2333 | 1.3887 | 1.5673 | **1.7623** | 1.6412 | 1.6870 |
| | | 0.5 | 1.0135 | 1.0824 | 1.1980 | 1.3494 | 1.5268 | **1.7223** | 1.5938 | 1.6394 |
| | | 1 | 1.0079 | 1.0713 | 1.1823 | 1.3301 | 1.5046 | **1.6986** | 1.5785 | 1.6185 |
| | | *2* | *1.0062* | *1.0666* | *1.1741* | *1.3182* | *1.4890* | ***1.6794*** | *1.5222* | *1.6088* |
| | | 5 | 1.0110 | 1.0748 | 1.1836 | 1.3274 | 1.4969 | **1.6852** | 1.5853 | 1.6322 |
| | | 10 | 1.0169 | 1.0862 | 1.1997 | 1.3471 | 1.5194 | **1.7098** | 1.6021 | 1.6465 |
| 1 | 20 | 0 | 1.0308 | 1.1164 | 1.2463 | 1.4091 | 1.5955 | **1.7992** | 1.9698 | ** |
| | | 0.5 | 1.0141 | 1.0880 | 1.2118 | 1.3732 | 1.5613 | **1.7688** | 1.6742 | 1.7282 |
| | | 1 | 1.0079 | 1.0765 | 1.1965 | 1.3554 | 1.5422 | **1.7489** | 1.6569 | 1.7090 |
| | | *2* | *1.0059* | *1.0721* | *1.1896* | *1.3461* | *1.5307* | ***1.7348*** | *1.6550* | *1.7052* |
| | | 5 | 1.0118 | 1.0830 | 1.2041 | 1.3629 | 1.5486 | **1.7539** | 1.6806 | 1.8899 |
| | | 10 | 1.0190 | 1.0979 | 1.2261 | 1.3910 | 1.5820 | 1.7911 | **2.0096** | 1.7653 |
| 2 | 10 | 0 | 1.0380 | 1.1439 | **1.3057** | 1.2902 | 1.3203 | 1.3534 | 1.3848 | 1.4118 |
| | | 0.5 | 1.0192 | 1.1105 | 1.2625 | **1.4619** | 1.2903 | 1.3230 | 1.3536 | 1.3799 |
| | | 1 | 1.0124 | 1.0970 | 1.2434 | **1.4380** | 1.2776 | 1.3099 | 1.3404 | 1.3656 |
| | | *2* | *1.0101* | *1.0912* | *1.2332* | ***1.4319*** | *1.2721* | *1.3040* | *1.3343* | *1.3592* |
| | | 5 | 1.0161 | 1.1007 | 1.2438 | **1.4329** | 1.2822 | 1.3136 | 1.3434 | 1.3695 |
| | | 10 | 1.0230 | 1.1140 | 1.2624 | **1.4662** | 1.2950 | 1.3265 | 1.3575 | 1.3836 |
| 2 | 20 | 0 | 1.0361 | 1.1366 | 1.2900 | **1.4841** | 1.3441 | 1.3777 | 1.4129 | 1.4454 |
| | | 0.5 | 1.0168 | 1.1027 | 1.2478 | **1.4374** | 1.3484 | 1.3484 | 1.3843 | 1.4176 |
| | | 1 | 1.0098 | 1.0891 | 1.2282 | **1.4141** | 1.3358 | 1.3358 | 1.3718 | 1.4056 |
| | | *2* | *1.0075* | *1.0833* | *1.2182* | ***1.3996*** | *1.3304* | *1.3304* | *1.3663* | *1.4003* |
| | | 5 | 1.0137 | 1.0937 | 1.2306 | **1.4120** | 1.3420 | 1.3420 | 1.3771 | 1.4118 |
| | | 10 | 1.0213 | 1.1085 | 1.2515 | **1.4377** | 1.3573 | 1.3573 | 1.3926 | 1.4269 |

**Table 5. Effect of aspect and thickness ratios on the nonlinear frequency ratio ($\omega_{NL}/\omega_L$) of simply supported FGM plates with temperature gradient ($T_c$ =600K, $T_m$ =300K)**

| a/b | a/h | k | *w/h* 0.2 | 0.4 | 0.6 | 0.8 | 1.0 | 1.2 | 1.4 | 1.6 |
|---|---|---|---|---|---|---|---|---|---|---|
| 1 | 10 | 0 | 1.0332 | 1.1218 | 1.2550 | 1.4212 | 1.6109 | **1.8176** | 1.6777 | 1.7246 |
| | | 0.5 | 1.0169 | 1.0941 | 1.2212 | 1.3856 | 1.5765 | **1.7860** | 1.6374 | 1.6855 |
| | | 1 | 1.0110 | 1.0830 | 1.2062 | 1.3680 | 1.5573 | **1.7661** | *1.6211* | *1.6706* |
| | | 2 | *1.0092* | *1.0789* | *1.1996* | *1.3588* | *1.5456* | ***1.7522*** | 1.6215 | 1.6710 |
| | | 5 | 1.0153 | 1.0899 | 1.2139 | 1.3750 | 1.5627 | **1.7693** | 1.6415 | 1.6888 |
| | | 10 | 1.0226 | 1.1047 | 1.2354 | 1.4023 | 1.5947 | **1.8055** | 1.6750 | 1.7194 |
| 2 | 20 | 0 | 1.0434 | 1.1586 | 1.3311 | **1.5462** | 1.3775 | 1.4146 | 1.4532 | 1.4788 |
| | | 0.5 | 1.0224 | 1.1243 | 1.2914 | **1.5072** | 1.3511 | 1.3767 | 1.4299 | 1.4666 |
| | | 1 | 1.0145 | 1.1104 | 1.2737 | **1.4876** | 1.3408 | 1.3595 | 1.4208 | 1.4576 |
| | | *2* | *1.0122* | *1.1058* | *1.2675* | ***1.4792*** | *1.3401* | *1.3520* | *1.4198* | *1.4583* |
| | | 5 | 1.0208 | 1.1224 | 1.2901 | **1.5070** | 1.3616 | 1.4011 | 1.4422 | 1.4801 |
| | | 10 | 1.0312 | 1.1444 | 1.3230 | **1.5495** | 1.3859 | 1.4264 | 1.4680 | 1.5063 |

**Table 6. Influence of aspect ratio on the nonlinear frequency ratio ($\omega_{NL}/\omega_L$) of clamped FGM plates (a/h=20) with temperature gradient ($T_c$ =400K, $T_m$ =300K)**

| a/b | k | *w/h* | | | | |
|---|---|---|---|---|---|---|
| | | 0.2 | 0.4 | 0.6 | 0.8 | 1.0 |
| 1 | 0 | 1.0102 | 1.0403 | 1.0884 | 1.1521 | 1.2291 |
| | 0.5 | 1.0106 | 1.0417 | 1.0915 | 1.1572 | 1.2366 |
| | 1 | 1.0106 | 1.0417 | 1.0913 | 1.1570 | 1.2362 |
| | 2 | 1.0103 | 1.0408 | 1.0895 | 1.1540 | 1.2318 |
| | *5* | *1.0101* | *1.0398* | *1.0874* | *1.1505* | *1.2264* |
| | 10 | 1.0102 | 1.0401 | 1.0880 | 1.1513 | 1.2277 |
| 2 | 0 | 1.0109 | 1.0432 | 1.0954 | 1.1654 | 1.2509 |
| | 0.5 | 1.0112 | 1.0443 | 1.0978 | 1.1696 | 1.2570 |
| | 1 | 1.0111 | 1.0441 | 1.0973 | 1.1687 | 1.2557 |
| | 2 | 1.0109 | 1.0430 | 1.0950 | 1.1649 | 1.2501 |
| | *5* | *1.0106* | *1.0418* | *1.0924* | *1.1604* | *1.2435* |
| | 10 | 1.0106 | 1.0419 | 1.0926 | 1.1607 | 1.2439 |

**Table 7. Influence of skew angle of simply supported FGM plate (a/h=10, a/b=1) with temperature gradient (Tc =400K, Tm =300K) on the nonlinear frequency ratio ($\omega_{NL}/\omega_L$) .**

| Skew angle | k | *w/h* 0.2 | 0.6 | 1 | 1.2 | 1.6 |
|---|---|---|---|---|---|---|
| $15^0$ | 0 | 1.0291 | 1.2342 | 1.5694 | **1.7661** | 1.6506 |
| | 0.5 | 1.0136 | 1.1989 | 1.5293 | **1.7258** | 1.6076 |
| | 1 | 1.0080 | 1.1831 | 1.5065 | **1.7012** | 1.5880 |
| | *2* | *1.0062* | *1.1751* | *1.4909* | ***1.6821*** | *1.5790* |
| | 5 | 1.0112 | 1.1846 | 1.4987 | **1.6879** | 1.5932 |
| | 10 | 1.0170 | 1.2004 | 1.5210 | **1.7123** | 1.6143 |
| $30^0$ | 0 | 1.0293 | 1.2788 | 1.5727 | **1.7709** | 1.5929 |
| | 0.5 | 1.0139 | 1.2003 | 1.5322 | **1.7292** | 1.5524 |
| | 1 | 1.0084 | 1.1846 | 1.5099 | **1.7058** | 1.5340 |
| | *2* | *1.0066* | *1.1765* | *1.4945* | ***1.6863*** | *1.5253* |
| | 5 | 1.0114 | 1.1857 | 1.5018 | **1.6914** | 1.5314 |
| | 10 | 1.0172 | 1.2012 | 1.5231 | **1.7144** | 1.5573 |
| $45^0$ | 0 | 1.0291 | 1.2354 | 1.5743 | **1.7676** | 1.4928 |
| | 0.5 | 1.0145 | 1.2017 | **1.5366** | 1.3869 | 2.1117 |
| | 1 | 1.0092 | 1.1868 | **1.5143** | 1.3731 | 1.4405 |
| | *2* | *1.0075* | *1.1787* | ***1.4993*** | *1.3659* | 1.4330 |
| | 5 | 1.0120 | 1.1871 | **1.5043** | 1.3771 | 1.4440 |
| | 10 | 1.0175 | 1.2017 | **1.5243** | 1.3927 | 1.4606 |

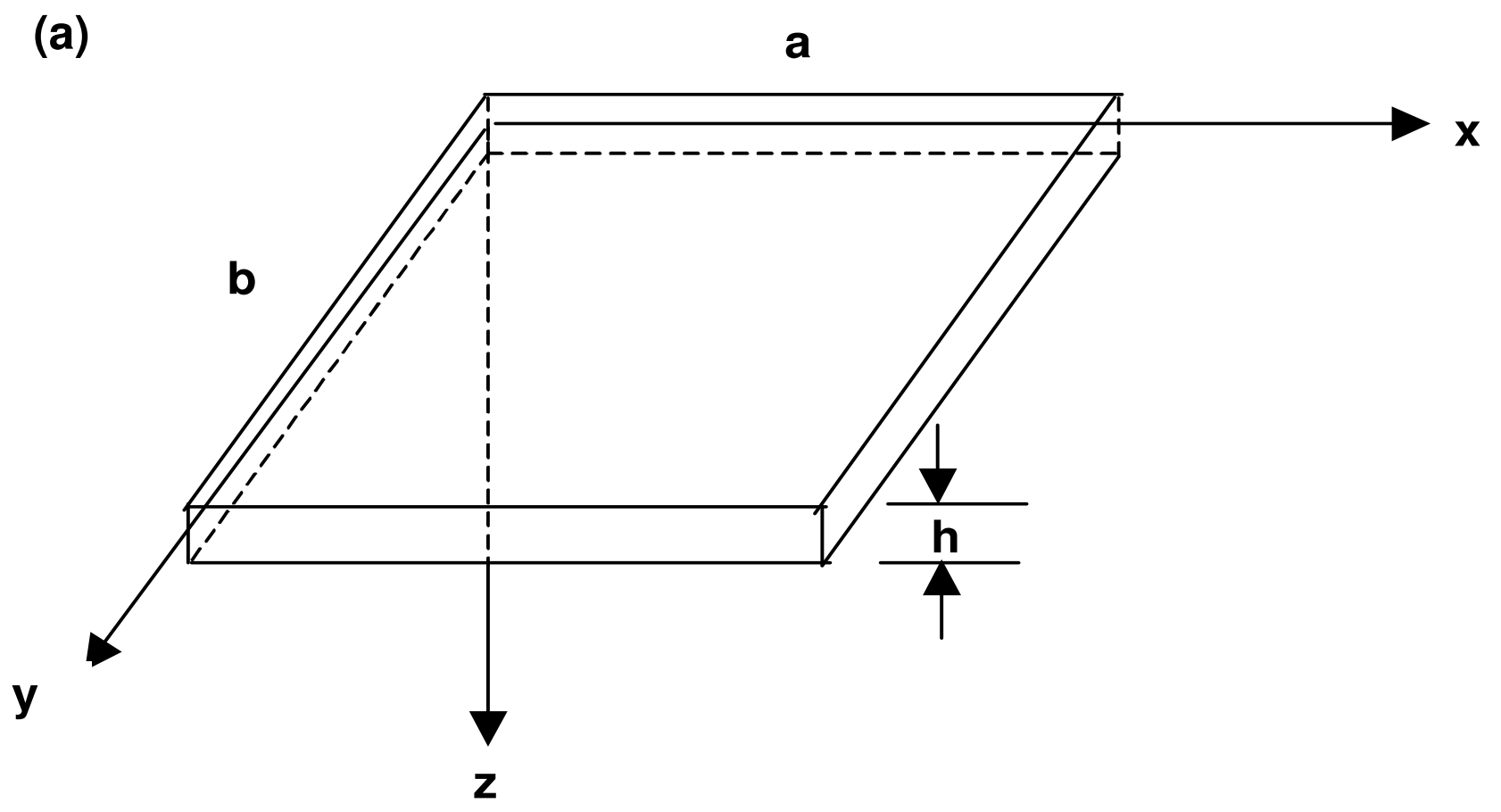

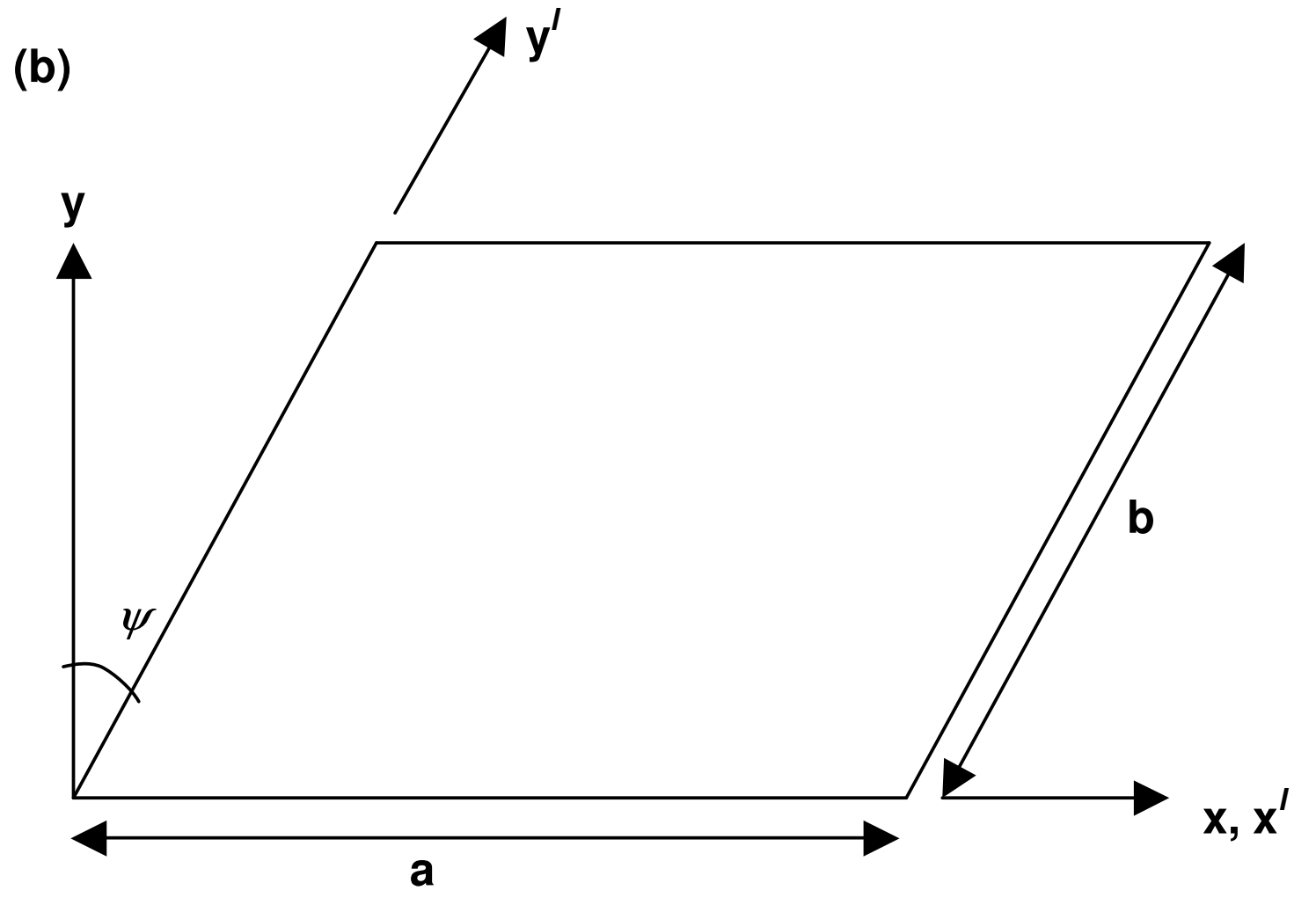

**Fig 1 (a). Configuration and coordinate system of a rectangular FGM plate.**
**(b). Coordinate system of a skew plate.**

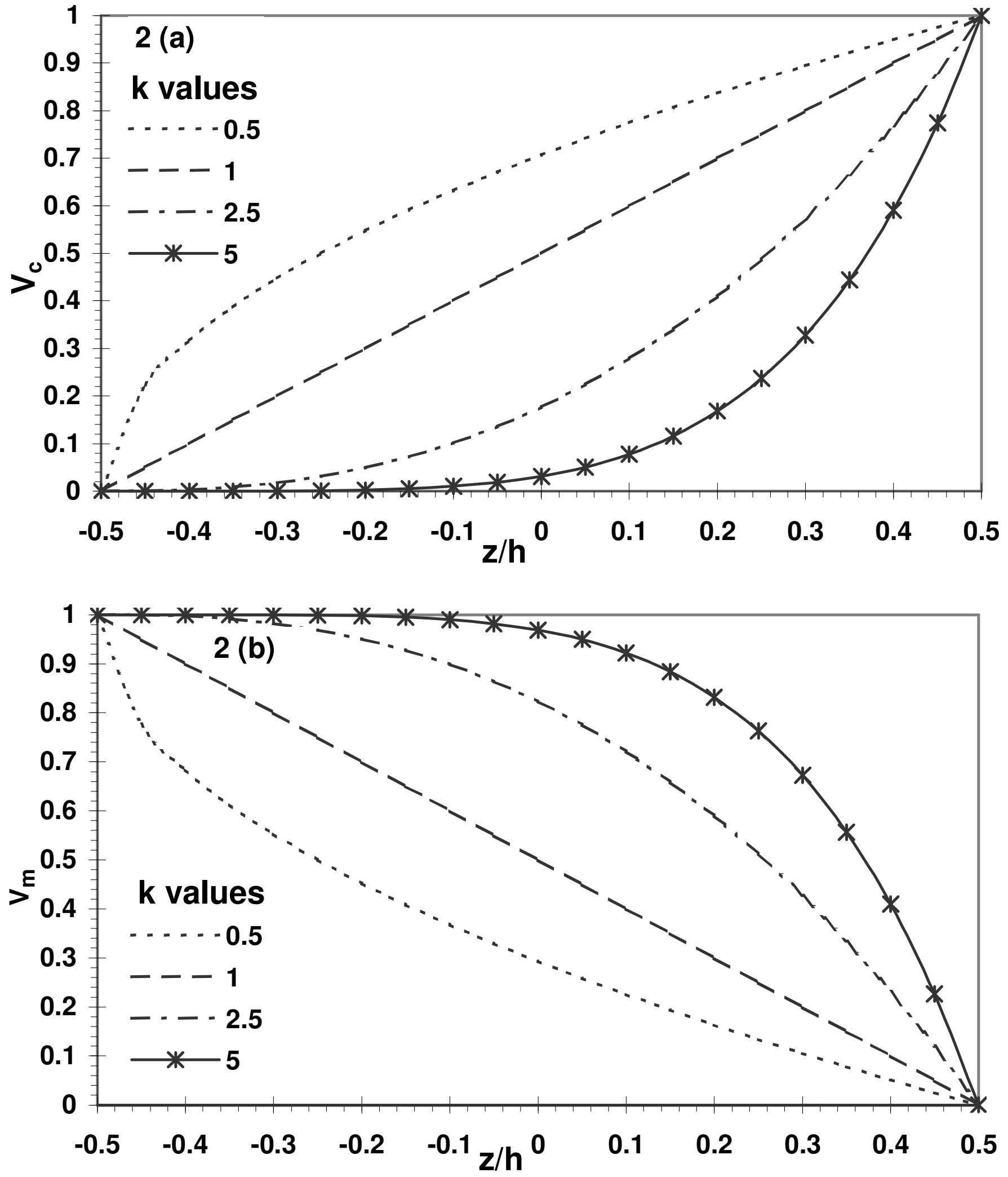

**Fig 2. Variation of volume fraction through the thickness: a) Ceramic; b) Metal**

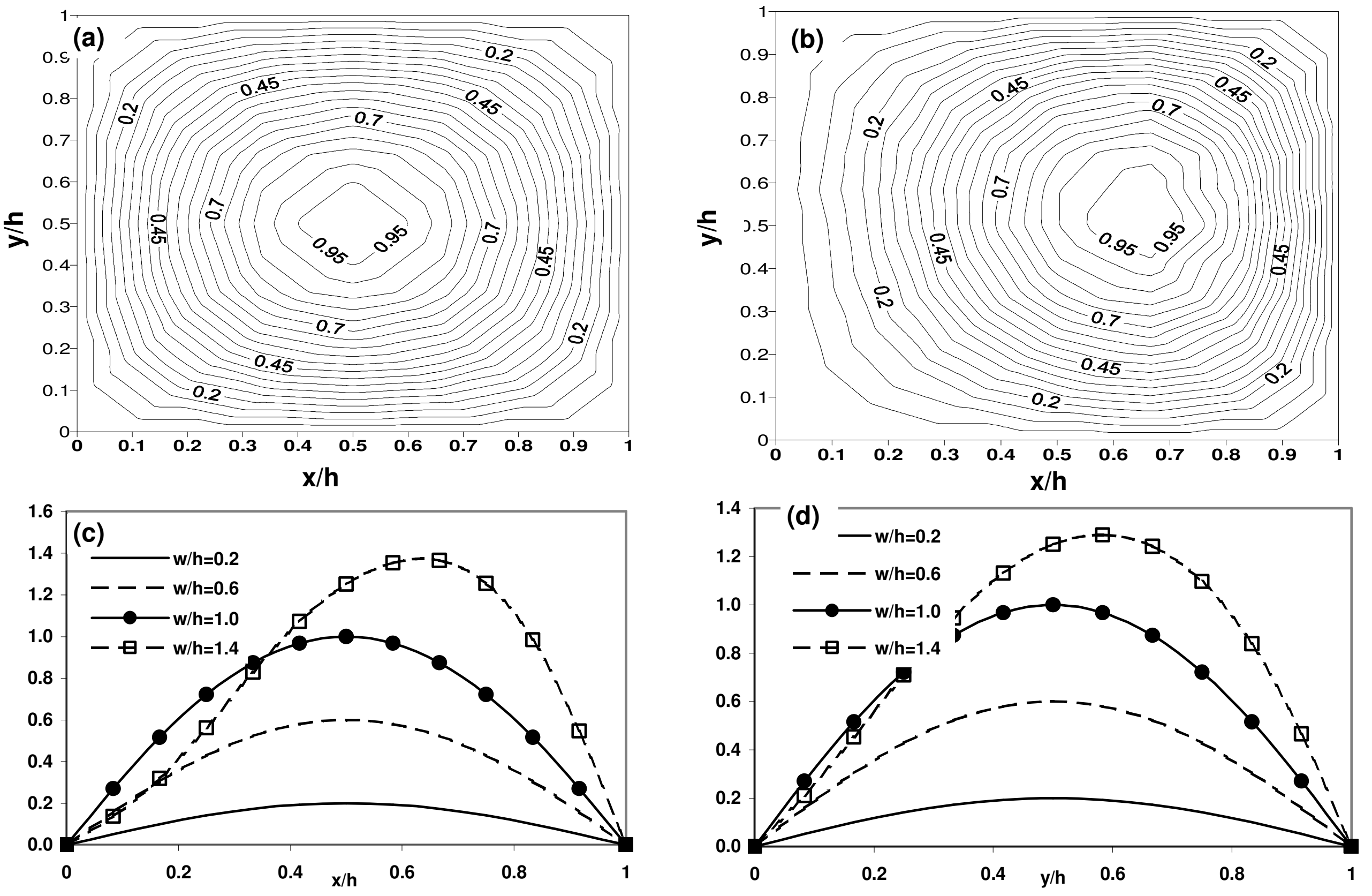

**Fig 3. The redistribution of normalized nonlinear mode shape contours of simply supported FGM plate (a/b=1, a/h=10, k=2): (a) w/h= 1.0; (b) w/h= 1.4; (c) mode shapes along y = b/2; (d) mode shapes along x = a/2.**